\documentclass[12pt,reqno,a4wide]{article}
\usepackage{makeidx}
\usepackage{comment}
\usepackage{mathptmx}
\usepackage{a4wide}
\usepackage{amssymb}
\usepackage{amsfonts}
\usepackage{amsthm}
\usepackage{amsmath}
\usepackage{dsfont}
\usepackage{graphicx}
\usepackage{shadow}
\usepackage[all]{xy}
\usepackage{enumerate}
\usepackage{mathrsfs}
\usepackage{upgreek}
\usepackage{stmaryrd}
\usepackage[utf8]{inputenc}
\usepackage{array}
\usepackage{tikz-cd}
\usepackage{thumbpdf}
\usepackage[colorlinks=true,pagebackref]{hyperref}
\usepackage{textcomp}
\usepackage{chapterbib}
\usepackage{graphics}
\usepackage{longtable}
\usepackage{epsf}
\usepackage{bm}
\usepackage{adjustbox}

\newtheorem{theorem}{Theorem}[section]
\newtheorem{lem}[theorem]{Lemma}
\newtheorem{thm}[theorem]{Theorem}
\newtheorem{prop}[theorem]{Proposition}
\newtheorem{cor}[theorem]{Corollary}
\theoremstyle{definition}
\newtheorem*{Beweis}{Proof}
\newtheorem{defn}[theorem]{Definition}
\newtheorem{ex}[theorem]{Example}
\newtheorem{exs}[theorem]{Examples}

\newtheorem{rem}[theorem]{Remark}
\newtheorem{rems}[theorem]{Remarks}
\newtheorem{punto}[theorem]{}

\input{tcilatex}
\swapnumbers
\CompileMatrices
\input{tcilatex}
\begin{document}

\title{Second Representable Modules over Commutative Rings
}
\author{ $%
\begin{array}{ccc}
\text{Jawad Abuhlail}\thanks{%
Corresponding Author; Email: abuhlail@kfupm.edu.sa. \newline
The authors would like to acknowledge the support provided by the Deanship
of Scientific Research (DSR) at King Fahd University of Petroleum $\&$
Minerals (KFUPM) for funding this work through projects No. RG1213-1 $\&$
RG1213-2} &  & \text{Hamza Hroub}\thanks{%
The paper is extracted from the Ph.D. dissertation of Dr. Hamza Hroub under
the supervision of Prof. Jawad Abuhlail} \\
\text{Department of Mathematics and Statistics} &  & \text{Department of
Mathematics } \\
\text{King Fahd University of Petroleum }\&\text{ Minerals} &  & \text{King
Saud University} \\
\text{31261 Dhahran, KSA} &  & \text{11451 Riyad, KSA}%
\end{array}%
$}
\date{\today }
\maketitle

\begin{abstract}
Let $R$ be a commutative ring. We investigate $R$-modules which can be written as \emph{finite} sums of {\it {second}} $R$-submodules (we call them \emph{second representable}). We provide sufficient conditions for an $R$-module $M$ to be have a (minimal) second presentation, in particular within the class of lifting modules. Moreover, we investigate the class of (\emph{main}) \emph{second attached prime ideals} related to a module with such a presentation.
\end{abstract}

\section*{Introduction}

    Throughout, $R$ is a commutative ring. We consider \emph{second representable} modules,
\emph{i.e. }modules which can be written as finite sums $M=\sum\limits_{i=1}^{n}M_{i}$ of second
submodules $M_{1},\cdots ,M_{n}$ of $_{R}M$ (recall that $N\leq M$ is said
to be \emph{second }iff $IN=N$ or $IN=0$ for every ideal $I\leq R$ \cite{Y2001}, \cite{Abu}).
The paper is divided in three sections. In Section 1, and for the convenience of the reader, we collect some preliminaries from Module Theory. In Section 2, is devoted to the study of second representable modules. In particular, we provide sufficient conditions for the existence of second representations for $_{R}M$ are
provided, among others, in Proposition \ref{Proposition 4.9} and Theorem \ref{Theorem 4.18}.

Second and semisimple modules are trivially second representable, and Example \ref{sr-not-ss} provides examples modules which are second representable but neither second not semisimple. Several other examples showing that some of the sufficient conditions in
the results mentioned above are not necessary (\emph{e.g.} Examples \ref{AP}
and \ref{Zn[x]}). Since every second module is secondary, the First and the Second Uniqueness
Theorems (Theorems \ref{Remark 4.4} and \ref{Remark 4.5}, respectively) for
a second representable $R$-module follow from the corresponding ones on
secondary representations (\cite{K1973}, \cite{B2009}). As a byproduct, we introduced a new class of modules lying properly between the classes of semisimple and lifting modules, namely the class of \emph{s-lifting modules} (see Figure 1 at the end of Section 2).

Section 3 is devote to the study of \emph{second attached prime ideals} for a given second representable $R$-module. As a consequence of Theorem \ref{Theorem 4.23}, it follows that a second representable Noetherian $R$-module is a finite \emph{direct} sum of second submodules. Theorem \ref{Theorem 4.25} investigates the relation between semisimple, multiplication and second representable modules.

\section{Preliminaries}
In this section, and for the convenience of the reader, we collect some definitions and results from the literature.

Throughout this paper, $R$ is a commutative ring, $M$ a non-zero $R$-module, $LAT(_{R}M)$ is the canonical lattice of
$R$-submodules of $M$ and we write $N \leq M$ to indicate that $N \in LAT(_{R}M)$.
For $N,K\leq M$ and $I\leq R,$ we set%
\begin{equation*}
(K:_{R}N):=\{r\in R\mid rN\subseteq K\}\text{ and }(N:_{M}I):=\{x\in M\mid
Ix\subseteq N\}.
\end{equation*}%
In particular, we set $Ann(N):=(0:_{R}N).$

\begin{punto} (\cite{Y2001}, \cite{Abu}) Let $M$ be an $R$-module.
An $R$-submodule $K \leq_R M$ is said to be \emph{second} \cite{Abu}
iff $K\neq 0$ and for any ideal $I\leq R$ we have
\begin{equation*}
IK=K\text{ or }IK=0.
\end{equation*}%
By $Spec^{s}(M)$, we denote the spectrum of second $R$-submodules of $M.$
\end{punto}

\begin{punto}
(\cite[Sec. 41]{W1991})\ We say that an $R$-submodule $N\leq M$ has a \emph{%
supplement} $K$ in $M$ iff there is an $R$-submodule $K\leq M$ minimal with
respect to $N+K=M$. The $R$-module $M$ is said to be \emph{supplemented} iff
every $R$-submodule of $M$ has a supplement in $M$. We say that $N\leq M$
has \emph{ample supplements} in $M$ \cite{W1991} iff for each submodule $%
U\leq M$ with $N+U=M$ there is a supplement $K\subseteq U$ of $N$ in $M$.
The $R$-module $M$ is called\emph{\ amply supplemented} iff every $R$%
-submodule of $M$ has ample supplements in $M$. For example, every Artinian
module is amply supplemented.
\end{punto}

\begin{punto}
A submodule $N\leq M$ is called \emph{small} (or \emph{superfluous}) \emph{%
in }$M$ \cite[19.1]{W1991} iff $N+K\neq M$ for any $R$-submodule $K\lneq M$.
An epimorphism of $R$-modules $f:M\longrightarrow M^{\prime }$ is called a
\emph{small epimorphism} iff $Ker(f)$ is small $M.$ An $R$-submodule $N\leq
M $ is called \emph{large} (or \emph{essential}\textit{)} \cite[17.1]{W1991}
iff $N\cap K\neq 0$ for any $R$-submodule $0\neq K\leq M$. A monomorphism of
$R$-modules $g:M\longrightarrow M^{\prime }$ is called a \emph{large
monomorphism} iff $f(M)$ is large in $M^{\prime }.$
\end{punto}

\begin{punto}
We say that $M$ is a \emph{lifting} $R$-module \cite[22.2]{JCNR} iff any $R$%
-submodule $N\leq M$ contains a direct summand $X\leq M$ such that $N/X$ is
small in $M/X$. An $R$-module $M$ is called \emph{extending} \cite[p. 265]%
{JCNR} iff every nonzero submodule of $M$ is essential in a direct summand
of $M$.
\end{punto}

\begin{punto}
An $R$-module $M$ is called \emph{uniform} \cite{JCNR} iff every nonzero $R$%
-submodule of $M$ is large in $M$ (equivalently, $0\in LAT(_{R}M)$ is
irreducible). An $R$-module $M$ has \textit{uniform dimension} $n$ \cite%
{JCNR}, and we write $u.dim(M)=n,$ iff there exists a large monomorphism
from a direct sum of $n$ uniform $R$-modules to $M$. An $R$-module $M$ is
\emph{hollow} iff every proper $R$-submodule of $M$ is small in $M$
(equivalently, $1\in LAT(_{R}M)$ is hollow). We say that $M$ has \textit{%
hollow }\emph{dimension }$n$ \cite{JCNR} iff there exists a small
epimorphism from $M$ to a direct sum of $n$ hollow $R$-modules, in this case
we write $h.dim(M)=n$.
\end{punto}

\begin{lem}
\label{Theorem 1.21}(\cite[Proposition 22.11]{JCNR}) Let $M$ be a nonzero $R$%
-module with finite hollow (uniform) dimension.

\begin{enumerate}
\item If $_{R}M$ is lifting, then $M=\bigoplus\limits_{i=1}^{n}H_{i}$ where
each $H_{i}$ is a hollow $R$-module and $n=h.dim(M)$.

\item If $_{R}M$ is extending, then $M=\bigoplus\limits_{i=1}^{n}U_{i}$
where each $U_{i}$ is a uniform $R$-module and $n=u.dim(M)$.
\end{enumerate}
\end{lem}

\begin{lem}
\label{Theorem 1.24}(\cite[22.2]{JCNR}, \cite[20.34]{JCNR})

\begin{enumerate}
\item Every lifting $R$-module is amply supplemented.

\item The following are equivalent for an amply supplemented $R$-module $M:$

\begin{enumerate}
\item $M$ has finite hollow dimension.

\item $M$ has the DCC on supplements.

\item $M$ has the ACC on supplements.
\end{enumerate}
\end{enumerate}
\end{lem}

\begin{punto}
\index{atomic}
\index{coatomic}
\label{min}Let $\mathcal{L}=(L,\wedge ,\vee, 0, 1)$ be a bounded lattice.

Let $X\subseteq L\backslash \{0\}$, and denote by $Min(X)$ the set of minimal elements of $X$.
We say that $X$ is \emph{atomic} iff for every $p\in X$ there is $q\in Min(X)$ such that $q\leq
p;$

Let $X\subseteq L\backslash \{1\}$, and denote by $Max(X)$ the set of maximal elements $X$.
We say that $X$ is \emph{coatomic} iff for every element $p\in X$ there is $q\in Max(X)$ such
that $p\leq q$.

Let $M$ be an $R$-module. We say that $M$ is \emph{atomic} (resp. \emph{coatomic}) iff the class of non-zero (resp. proper $R$ submodules) of $M$ is atomic (resp. coatomic).
\end{punto}

\begin{lem}
\label{Theorem 1.26}(\cite[41.5, 41.6]{W1991})

\begin{enumerate}
\item If $_{R}M$ is coatomic and every maximal $R$-submodule of $M$ has a
supplement in $M,$ then $M$ is a sum of hollow submodules of $M$.

\item Let $_{R}M$ be finitely generated. Then $M$ is supplemented if and
only if $M$ is a sum of hollow submodules.
\end{enumerate}
\end{lem}

\subsection*{Primary and Secondary Representations}

\begin{punto}
\index{primary submodule}
\index{primary decomposition}
\index{minimal primary submodule} A proper $R$-submodule $N\lneq M$ is
called \emph{primary} \cite{AK2012} iff whenever $ax\in N$ and $x\notin N$
we have $a^{n}M\subseteq N$ for some $n\in \mathbb{N}$. If $N$ is a primary
submodule of $M,$ then $p:=%
\sqrt{(N:_{R}M)}$ is prime ideal of $R$ and we say that $N$ is $p$\emph{%
-primary}. A submodule $K\leq M$ has a \emph{primary decomposition} \cite%
{AK2012} iff there are primary submodules $N_{1},\cdots ,N_{n}$ of $M$ with $%
K=\bigcap_{i=1}^{n}N_{i}.$ Such a decomposition of $K,$ if it exists, is
called \emph{minimal} iff:

\begin{enumerate}
\item $\sqrt{(N_{i}:_{R}M)}\neq \sqrt{(N_{j}:_{R}M)}$ for $i\neq j;$

\item $\bigcap\limits_{i\neq j}N_{i}\nsubseteq N_{j}$ $\forall j\in
\{1,2,\cdots ,n\}.$
\end{enumerate}
\end{punto}

\begin{thm}
\label{Theorem 1.17}(Lasker-Noether Theorem \cite[Theorem 18.20]{AK2012})
Every submodule of a finitely generated module over a Noetherian ring has a
primary decomposition.
\end{thm}

\begin{thm}
\label{Theorem 1.15}(First uniqueness Theorem of Primary Decompositions \cite%
[Theorem 18.19]{AK2012}) Let $R$ be Noetherian and $M$ an $R$-module. If $%
\bigcap\limits_{i=1}^{n}N_{i}=N=\bigcap\limits_{j=1}^{m}K_{j}$ are two
minimal primary decompositions of $N\leq M$, where $N_{i}$ is $p_{i}$-
primary for all $i\in \{1,2,\cdots ,n\}$ and $K_{j}$ is $q_{j}$-primary for
all $j\in \{1,\cdots ,m\}$, then $n=m$ and $\{p_{1}\cdots
,p_{n}\}=\{q_{1},\cdots ,q_{m}\}$.
\end{thm}

\begin{thm}
\label{Theorem 1.16}(Second Uniqueness Theorem of Primary Decomposition (%
\cite[Theorem 18.24]{AK2012}) Let $M$ be a finitely generated module over a
Noetherian ring $R$ and $\bigcap\limits_{i=1}^{n}N_{i}=N=\bigcap%
\limits_{i=1}^{n}K_{i}$ be two minimal primary decompositions of $N\leq M$,
where $N_{i}$ and $K_{i}$ are $p_{i}$- primary submodules of $M$ for all $%
i\in \{1,\cdots ,n\}$. If $p_{j}$ is minimal among $\{p_{1},\cdots
,p_{n}\}$ for some $j\in \{1,\cdots ,n\},$ then $N_{j}=K_{j}$.
\end{thm}

Dual to primary submodules and primary decompositions are secondary submodules and secondary representations:

\begin{punto}
\index{representable module}
\index{attached prime} An $R$-module $M$ is called \emph{secondary}%
\index{secondary module} (\cite{K1973}, \cite{M1973}) iff for any $a\in R$
we have $aM=M$ or $a^{n}M=0$ for some $n\in \mathbb{N}$. If $M$ is a
secondary $R$-module, then $p:=%
\sqrt{Ann(M)}$ is a prime ideal of $R$ and $M$ is called $p$\emph{-secondary}%
. An $R$-module $M$ is called \emph{representable} (\cite{K1973}, \cite%
{M1973}) iff $M=\sum\limits_{i=1}^{n}N_{i},$ where $N_{1},\cdots ,N_{n}$ are
secondary $R$-module. Moreover, $M=\sum\limits_{i=1}^{n}N_{i}$ is said to
be a \emph{minimal secondary representation}
\index{minimal secondary representation} iff $%
\sqrt{Ann(N_{i})}\neq \sqrt{Ann(N_{j})}$ whenever $i\neq j$ and $%
N_{j}\nsubseteq \sum\limits_{i\neq j}N_{i}$ for all $j\in \{1,\cdots ,n\}$.
For each $i\in \{1,\cdots ,n\},$ the prime ideal $p_{i}:=\sqrt{Ann(N_{i})}$
is called an \emph{attached prime} \cite{M1973} and we set $%
Att(M):=\{p_{1},\cdots ,p_{n}\}.$ A subset $A\subseteq Att(M)$ is called
\emph{isolated} iff $q\in A$ whenever $q\in Att(M)$ and $q\subseteq p$ for
some $p\in A$. Examples of representable modules are artinian modules (\cite[Theorem 1]{K1973})
and injective modules over Noetherian rings (\cite[Theorem 2.3]{S1976}). Every quotient $Q$ of a
representable module $_{R}M$ is representable and $Att(Q)\subseteq Att(M)$ (\cite[Theorem 1.10]{Y1995}).
\end{punto}

\begin{punto}
\index{coassociated prime} A prime ideal $p\leq R$ is called a \emph{%
coassociated prime} \cite{C1981} to $_{R}M$ iff there is a \emph{hollow}
factor $M^{\prime }$ of $M$ such that $p=\{a\in R\mid aM^{\prime }\neq
M^{\prime }\}$. The set of coassociated primes of an $R$-module $M$ is
denoted by $Coass(M)$. If $_{R}M$ is representable, then $Att(M)=Coass(M)$ (%
\cite[Theorem 1.14]{Y1995}).
\end{punto}

\begin{prop}
\label{Proposition 1.29}(\cite{M1973}) Let $M_{1},\cdots ,M_{n}$ be
secondary $R$-submodules of the $R$-module $M.$ Then $M_{1}\oplus
M_{2}\oplus \cdots \oplus M_{n}$ is a $p$-secondary $R$-module if and only
if $M_{i}$ is a $p$-secondary $R$-submodule of $M$ for all $i$.
\end{prop}

\begin{thm}
\label{Theorem 1.30}(First Uniqueness Theorem of Secondary Representations) (%
\cite[Theorem 2]{K1973}) If $\sum\limits_{i=1}^{n}K_{i}$ and $%
\sum\limits_{j=1}^{m}N_{j}$ are two minimal secondary representations for $%
_{R}M,$ with $K_{i}$ is $p_{i}$ -secondary for $i=1,\cdots ,n$ and $N_{j}$
is $q_{j}$-secondary for $j=1,\cdots ,m$, then $n=m$ and $%
\{p_{1},\cdots ,p_{n}\}=\{q_{1},\cdots ,q_{n}\}$.
\end{thm}

\begin{thm}
\label{Theorem 1.31} (Second Uniqueness Theorem of Secondary Representation
\cite[Theorem 3.2.7]{B2009}) Let $M$ be representable, $A\subseteq Att(M)$
an isolated subset and $M=\sum\limits_{i=1}^{n}K_{i}$ a minimal secondary
representation for $M$ with $K_{i}$ is $p_{i}$-secondary for $i\in
\{1,\cdots ,n\}.$ Then $\sum\limits_{p_{i}\in A}K_{i}$ is independent of the
choice of the minimal second representation.
\end{thm}

\begin{thm}
\label{Theorem 1.33} (\cite[Theorem 1.10]{Y1995}) If $N$ is a representable $%
R$-submodule of the representable module $M$, then $Att(M/N)\subseteq
Att(M)\subseteq Att(N)\cup Att(M/N)$.
\end{thm}

\begin{thm}
\label{Theorem 1.34}(\cite[Theorem 1.11]{Y1995}) If $M_{1},\cdots ,M_{n}$
are representable $R$-modules, then $\bigoplus\limits_{i=1}^{n}M_{i}$ is
representable and $Att(\bigoplus\limits_{i=1}^{n}M_{i})=\bigcup%
\limits_{i=1}^{n}Att(M_{i})$.
\end{thm}

\section{Second Representations}

Recall that $R$ is a commutative ring. Yassemi \cite{Y2001} introduced the notion of \emph{second submodules} of a
given non-zero module over a commutative ring. Annin \cite{Ann2002} called
these \emph{coprime modules} (see also \cite{Wij2006}) and used them to
dualize the notion of attached primes.

\begin{punto}
\label{second}A nonzero submodule $K\leq M$ is called \emph{second}%
\index{second submodule} \cite{Y2001} iff for any ideal $I\leq R$, we have $%
IK=0$ or $IK=K $. The spectrum of second $R$-submodules of $M$ is denoted by
$Spec^{s}(M).$ If $K\in Spec^{s}(M),$ then $p:=(0:_{R}K)$ is a prime ideal,
called a \emph{second attached prime}%
\index{second attached prime}\emph{\ }of $M$ and $K$ is called $p$\emph{%
-second}%
\index{p-second submodule}. By%
\begin{equation*}
Att^{s}(M):=\{(0:_{R}K)\mid K\in Spec^{s}(M)\}  \label{Att^s(M)}
\end{equation*}%
we denote the set of second attached primes of $M.$
\end{punto}

\begin{lem}
\label{Proposition 4.2} Let $\{K_{i}\}_{i\in A}$ be family of second $R$%
-submodule of $M$ such that $K_{j}\nsubseteq \sum\limits_{i\in A\backslash
\{j\}}K_{i}$ for all $j\in A$. Let $p$ be a prime ideal of $R.$ Then $K_{i}$
is $p$-second for all $i\in A$ if and only if $\sum_{i\in A}K_{i}$ is $p$%
-second.
\end{lem}

\begin{Beweis}
$(\Rightarrow )$ Assume that $K_{i}$ is $p$-second for all $i\in A$.
Clearly, $p=(0:_{R}\sum\limits_{i\in A}K_{i})$. Let $I\leq R.$ If $IK_{j}=0$
for some $j\in A,$ then $I\subseteq p,$ whence $I\sum_{i\in A}K_{i}=0$.
Otherwise, $IK_{j}=K_{j}$ for all $j\in A$ and so $I\sum\limits_{i\in
A}K_{i}=\sum\limits_{i\in A}IK_{i}=\sum\limits_{i\in A}K_{i}$. Consequently,
$\sum_{i\in A}K_{i}$ is second.

$(\Leftarrow )$ Assume that $\sum\limits_{i\in A}K_{i}$ is $p$-second and
that $K_{i}$ is $p_{i}$-second for $i\in A$. Clearly $p\subseteq p_{i}$ for
all $i\in A$. For any $j\in A$, we have $p_{j}\sum\limits_{i\in
A}K_{i}=\sum\limits_{i\in A\backslash \{j\}}p_{j}K_{i}\subseteq
\sum\limits_{i\in A\backslash \{j\}}K_{i}\neq \sum\limits_{i\in A}K_{i},$
whence $p_{j}\sum\limits_{i\in A}K_{i}=0$, i.e. $p_{j}\subseteq p$. Hence $%
p=p_{j}$ for all $j\in A$.$\blacksquare$
\end{Beweis}

\begin{defn}
\label{second rep}We say that an $R$-module $M$ is (\emph{directly}%
\index{directly second representable}) \emph{second representable}%
\index{second representable} iff $M=\sum\limits_{i=1}^{n}K_{i}$ ($%
M=\bigoplus\limits_{i=1}^{n}K_{i}$) where $K_{i}$ is a second $R$-submodule
of $M$ for all $i=1,2,\cdots ,n;$ in this case we call $\sum%
\limits_{i=1}^{n}K_{i}$ ($\bigoplus\limits_{i=1}^{n}K_{i}$) a (\emph{direct}%
\index{direct second representation}) \emph{second representation}%
\index{second representation} of $M$. An $R$-module $M$ is called \emph{%
semisecond}%
\index{semisecond} iff $M$ is a (not necessarily finite) sum of second
submodules of $M$.
\end{defn}

\begin{ex}
\item Let $p$ be a prime number. Any divisible $p$-group is a semisecond $%
\mathbb{Z}$-module but not semisimple. This follows from the fact that every
divisible $p$-group is a direst sum of copies of the Pr\"{u}fer group $\mathbb{Z}%
(p^{\infty })$ which is a $0$-second $\mathbb{Z}$-module but not simple (see
\cite[p. 124]{W1991}, \cite[p. 96]{PRU} for the properties of the Pr\"{u}fer group).
\end{ex}

\begin{punto}
\label{minimal second rep}A (direct) second representation $%
M=\sum\limits_{i=1}^{n}K_{i}$ ($M=\bigoplus\limits_{i=1}^{n}K_{i}$) is
called a \emph{minimal }(\emph{direct})\emph{\ second representation}%
\index{minimal second representation} for $M$ iff it satisfies the following
conditions:

\begin{enumerate}
\item $(0:_{R}K_{i})\neq (0:_{R}K_{j})$ for $i\neq j$.

\item $K_{j}\nsubseteq \sum\limits_{i=1,i\neq j}^{n}K_{i}$ for all $%
j=1,2,\cdots ,n$.
\end{enumerate}
\end{punto}

\begin{punto}Let $_{R}M$ be second representable. It is clear that $M$ has a minimal
second representation say $\sum\limits_{i=1}^{n}K_{i}.$ Each $K_{i}$ in such
a minimal representation is called a \textit{main second submodule}%
\index{main second submodule} of $M$ and $(0:_{R}K_{i})$ is called a \emph{%
main second attached prime}%
\index{main second attached prime} of $M.$ So, the set of main second
attached primes is%
\begin{equation*}
att^{s}(M)=\{Ann(K_{i})\mid i=1,\cdots ,n\}.  \label{att^s(M)}
\end{equation*}%
By Theorem \ref{Remark 4.4} below, $att^{s}(M)$ is independent of the choice
of the minimal second representation of $M.$
\end{punto}

The result follows by Lemma \ref{Proposition 4.2}:

\begin{thm}
\label{existence of minimal second repersentation}(Existence Theorem for
Minimal Second Representations) Let $M $ be second representable. Then $M$
has a minimal second representation.
\end{thm}

\begin{ex}
\label{not-second-}The Abelian group $\mathbb{Z}_{18}$ has a minimal
secondary representation as a $\mathbb{Z}$-module, namely $\mathbb{Z}_{18}=2%
\mathbb{Z}_{18}+9\mathbb{Z}_{18}.$ However, $\mathbb{Z}_{18}$ has \emph{no}
second representation ($9\mathbb{Z}_{18}$ is the unique second $\mathbb{Z}$%
-submodule of $\mathbb{Z}_{18}$).
\end{ex}

In the light of \ref{second-seconday}, we obtain as a direct
consequence of Theorem \ref{Theorem 1.30} and \ref{Theorem 1.31} the First $%
\&\ $Second Uniqueness Theorems for Second Representations:

\begin{thm}
\label{Remark 4.4}(First Uniqueness Theorem of Second Representations) Let $%
M $ be an $R$-module with two minimal second representations $%
\sum\limits_{i=1}^{n}K_{i}=M=\sum_{j=1}^{m}N_{j},$ where $K_{i}$ is $p_{i}$%
-second for all $i\in \{1,2,\cdots ,n\}$ and $N_{j}$ is $q_{j}$-second for
all $j\in \{1,2,\cdots ,m\}$. Then $\{p_{1},p_{2},\cdots
,p_{n}\}=\{q_{1},q_{2},\cdots ,q_{m}\}$.
\end{thm}

\begin{thm}
\label{Remark 4.5}(Second Uniqueness Theorem of Second Representations) Let $%
M$ be second representable. If $\sum\limits_{i=1}^{n}K_{i}=M=\sum%
\limits_{i=1}^{n}N_{i}$ are minimal second representations for $M$ with $%
K_{j}$ and $N_{j}$ are $p_{j}$-second submodules of $M$ and $p_{j}$ is
minimal in $\{p_{1},\cdots ,p_{n}\}$ for some $j\in \{1,\cdots ,n\}$, then $%
K_{j}=N_{j}$.
\end{thm}

\begin{rems}
\index{second representable} \label{Remark 4.6}Let $M,$ $M_{1}, ..., M_{n}$
be second representable submodules of an $R$-module $L$.

\begin{enumerate}
\item $\sum\limits_{i=1}^{n}M_{i}$ is second representable and $%
att^{s}(\sum\limits_{i=1}^{n}M_{i})\subseteq
\bigcup\limits_{i=1}^{n}att^{s}(M_{i})$.

\item Any quotient $Q$ of $M$ is second representable and $%
att^{s}(Q)\subseteq att^{s}(M)$. To see this, let $K_{1}+\cdots +K_{n}$ be a minimal second representation for $M$
and $Q=M/N$ for some $R$-submodule $N\leq M$, then $M/N=\sum%
\limits_{i=1}^{n}(K_{i}+N)/N.$ It is easy to see that $(K_{i}+N)/N$ is
second for $i\in \{1,\cdots ,n\}$. The result is obtained now by applying
(1).

\item If $N$ is a second representable submodule of $M$, then it follows from (2) and Theorem \ref{Theorem 1.33} that
\begin{equation*}
att^{s}(M/N)\subseteq att^{s}(M)\subseteq att^{s}(N)\cup att^{s}(M/N).
\end{equation*}

\item If $M_{j}\cap \sum\limits_{i\neq j}M_{i}=0$ for all $j$, then it follows from Theorem \ref{Theorem 1.34} that
$\bigoplus_{i=1}^{n}M_{i}$ is second representable and
\begin{equation*}
att^{s}(\bigoplus\limits_{i=1}^{n}M_{i})=\bigcup%
\limits_{i=1}^{n}att^{s}(M_{i}).
\end{equation*}

\item For any multiplicatively closed subset of $S$ of $R$, the $S^{-1}R$-module $S^{-1}M$ is second representable and
\begin{equation*}
att^{s}(S^{-1}M)=\{p_{S}\mid p\in att^{s}(M)%
\text{ such that }p\cap S=\emptyset \}.
\end{equation*}
\end{enumerate}
\end{rems}

\begin{prop}
\index{second representable} \label{Remark 4.10}Let $M=\sum\limits_{i\in
\Lambda }K_{i}$ (resp. $M=\bigoplus\limits_{i\in \Lambda }K_{i}$), where $%
K_{i}$ is second for every $i\in \Lambda $. If $Att^{s}(M)$ is finite, then $%
M$ is second representable (resp. directly second representable).
\end{prop}

\begin{Beweis}
Let $M=\sum\limits_{i\in \Lambda }K_{i}$ (resp. $M=\bigoplus\limits_{i\in
\Lambda }K_{i}$) such that each $K_{i}$ is second for every $i\in \Lambda $.
Assume that $Att^{s}(M)=\{p_{1},p_{2},\cdots ,p_{n}\}$. For $j\in
\{1,2,\cdots ,n\}$, set $A_{j}=\{K_{i}:i\in \Lambda $ such that $%
(0:_{R}K_{i})=p_{j}\}.$ Notice that $N_{j}=\sum\limits_{K_{i}\in A_{j}}K_{i}$
is second by Proposition \ref{Proposition 4.2} for each $j\in \{1,2,\cdots
,n\}.$ Moreover, $M=\sum\limits_{j=1}^{n}N_{j}$ (resp. $M=\bigoplus%
\limits_{j=1}^{n}N_{j}$).$\blacksquare$
\end{Beweis}

\begin{punto}
We say that a submodule $K\leq M$ satisfies the
\index{IS-condition} \emph{IS-condition} iff for every $I\leq R$ for which $%
IK\neq 0,$ the submodule $IK\leq M$ has a proper supplement in $M.$
\end{punto}

\begin{rem}
\label{Remark 4.7}Let $_{R}M$ be supplemented, $K\leq M$ and $0\neq H\leq K$%
. The following conditions are equivalent:

\begin{enumerate}
\item $K$ is not contained in any supplement of $H$ in $M$.

\item $H$ has a proper supplement in $M$.
\end{enumerate}
\end{rem}

\begin{Beweis}
$(1\Rightarrow 2)$ Assume that $K$ is not contained in any supplement of $H$
in $M$. Since $M$ is supplemented, $H$ has a supplement $L$ in $M$, i.e. $%
H+L=M$. Indeed $L\neq M$ as $K\nsubseteq L$.

$(2\Rightarrow 1)$ Assume that $H$ has a proper supplement $L$ in $M$. Then $%
K\nsubseteq L;$ otherwise, $H+L=L\neq M$.$\blacksquare$
\end{Beweis}

\begin{lem}
\label{IS-Second}Every hollow $R$-submodule $0\neq K\leq M$ satisfying the
IS-condition is second.%
\end{lem}

\begin{Beweis}
Let $0\neq K\leq M$ be a hollow $R$-submodule satisfying the IS-condition.
Let $I\leq K$ and suppose that $IK\neq 0.$ By the IS-condition, $IK$ has
a proper supplement $L\lvertneqq M.$ It is easy to show that $IK+(L\cap K)=K$.
Since $K$ is hollow, $IK=K$ (notice that $L\cap K\neq K;$ otherwise, $%
IK+L=L\neq M$).$\blacksquare$
\end{Beweis}

\begin{ex}
\label{Example 4.8}The Abelian group $\mathbb{Z}_{18},$ considered as a $%
\mathbb{Z}$-module, is supplemented but not semisimple. The submodule $%
K_{1}=9\mathbb{Z}_{18}$ is hollow and satisfies the IS-condition. Notice
that $K_{2}:=6\mathbb{Z}_{18}$ is hollow and second but does not satisfy the
IS-condition (i.e. the IS-condition is not necessary for a hollow submodule
module to be second).
\end{ex}

\begin{punto}
\index{hollow representable}
\index{directly hollow representable} We say that an $R$-module $M$ is (%
\emph{directly}) \emph{hollow representable} iff $M$ is a finite (direct)
sum of hollow $R$-submodules.
\end{punto}

\begin{prop}
\index{directly hollow representable}
\index{hollow representable} \label{Proposition 4.9}Let $_{R}M$ be
(directly) hollow representable. If every \emph{maximal} hollow non-zero
submodule of $M$ is second, then $M$ is (directly) second representable.
\end{prop}

\begin{Beweis}
Let $M=\sum\limits_{i=1}^{n}H_{i}$ be a sum of hollow $R$-submodules.
Assume, without loss of generality, that this sum is \emph{irredundant.}
\textbf{Claim:} $H_{1},\cdots ,H_{n}$ are maximal hollow submodules of $M$.
To see this, suppose that $H$ is a hollow submodule of $M$ with $H_{i}\leq H$
for some $i\in \{1,\cdots ,n\}$ and consider $N:=\sum\limits_{j\neq i}H_{j}$%
. For any $x\in H$, there are $y\in N$ and $z\in H_{i}$ such that $x=y+z$.
But $z\in H$ implies that $y\in H$. So, $H=(N\cap H)+H_{i}$. Since $H$ is
hollow, either $H\cap N=H$ whence $H\subseteq N$, or $H_{i}=H.$ But $%
H_{i}\subseteq H$ and $M=\sum\limits_{i=1}^{n}H_{i}$ is an irredundant sum,
whence $H=H_{i}$. Hence $H_{i}$ is maximal hollow. By our assumption, $%
H_{1},\cdots ,H_{n}$ are second, whence $M=\sum\limits_{i=1}^{n}H_{i}$ is a
second representation of $M.$

If $M=\bigoplus\limits_{i=1}^{n}H_{i}$ is a direct sum of hollow $R$%
-submodules, then one can show similarly that each $H_{i}$ is a maximal
hollow $R$-submodule of $M,$ whence $M$ is a direct sum of hollow $R$-submodules.$\blacksquare$
\end{Beweis}

\begin{ex}
\label{Theorem 1.18} Every Artinian left $R$-module is hollow representable
(see \cite[Lemma 3.2]{ST2009}). Let $p\in \mathbb{Z}$ be a prime number. The
\emph{Pr\"{u}fer group}, considered as a $\mathbb{Z}$-module, is Artinian and
the unique maximal hollow $\mathbb{Z}$-submodule of $\mathbb{Z}({p^{\infty }}%
)$ is second.
\end{ex}

\begin{ex}
\label{lifting-hollow}A lifting $R$-module $M$ is directly hollow
representable if it satisfies any of the following additional conditions:

\begin{enumerate}
\item $_{R}M$ has a finite hollow dimension \cite[Proposition 22.11]{JCNR}
(e.g. $_{R}M$ is finitely generated \cite[Corollary 22.12]{JCNR}).

\item $_{R}M$ has a finite uniform dimension \cite[Proposition 22.11]{JCNR}
(e.g. $_{R}M$ is finitely cogenerated \cite[Corollary 22.12]{JCNR}).
\end{enumerate}
\end{ex}

Inspired by Example \ref{lifting-hollow} and Proposition \ref{Proposition
4.9}, we introduce the notion of s-lifting modules.

\begin{defn}
\label{stongly lifting}We call $_{R}M$ \emph{s-lifting}%
\index{s-lifting module} iff $_{R}M$ is lifting and every \emph{maximal}
hollow submodule of $M$ is second.
\end{defn}

\begin{exs}
\label{Example 4.14}

\begin{enumerate}
\item Consider the Abelian group $M=\mathbb{Z}_{8}$ as a $\mathbb{Z}$%
-module. Notice that $N=\{0,4\}$ is the unique second submodule in $M$,
hence $M$ is not second representable. Notice that $_{\mathbb{Z}}M$ is
Artinian and lifting but not s-lifting.

\item Every semisimple module is s-lifting and trivially semisecond (every
simple submodule is second).

\item Every second hollow module is s-lifting but not necessarily simple.
Consider the Pr\"{u}fer group $M=\mathbb{Z}({p^{\infty }})$, considered as $\mathbb{Z}$-module.
Notice that $_{\mathbb{Z}}M$ is not simple. Moreover, $_\mathbb{Z}M$ is
hollow and second whence s-lifting hollow but not semisimple.
\end{enumerate}
\end{exs}

As a direct consequence of Proposition \ref{Proposition 4.9} and Example \ref%
{lifting-hollow}, we obtain the following class of directly second
representable modules.

\begin{ex}
\index{s-lifting} \label{Theorem 4.15}If $_{R}M$ is an s-lifting module and
has a finite hollow dimension, then $M$ is directly second representable.
Clearly, this class is nonempty; \emph{e.g.} any finite direct sum of Pr\"{u}fer
groups is s-lifting with finite hollow dimension.
\end{ex}

The following example is an s-lifting second module with infinite hollow
dimension which is not semisimple.

\begin{ex}
\item \label{AP}Let $\mathbb{P}$ be the set of prime numbers, $A\subseteq
\mathbb{P}$ infinite and consider $M:=\bigoplus\limits_{p\in A}\mathbb{Z}%
(p^{\infty })$ considered as a $\mathbb{Z}$-module.

\textbf{Claim:}\ $_{\mathbb{Z}}M$ is lifting. This can be obtained by
applying \cite[Theorem 2]{BH1990} (see the second paragraph on page 60). However, we provide here our
own proof.

Let $N\leq M$. Assume, without loss of generality that $N$ is not a direct
summand of $M$ (if\textbf{\ }$N$ is a direct summand of $M$, then $N/N=0$ is
small in $M/N$ and we are done). Notice that $N=\bigoplus\limits_{p\in
A}L_{p},$ where $L_{p}\leq \mathbb{Z}(p^{\infty })$ for all $p\in A$.

\textbf{Case 1: }$L_{p}\neq Z(p^{\infty })$ for all $p\in A$. In this case, $%
N$ is small in $M$ as the set of submodules of $\mathbb{Z}(p^{\infty })$
form a chain for all $p\in A$. Indeed, for every $p\in A:$ if $L_{p}+W_{p}=M$%
, then $L_{p}\subseteq W_{p}=M$.

\textbf{Case 2: }$L_{p}=Z(p^{\infty })$ for all $p\in B\subsetneq A$ and $%
L_{p}\neq Z(p^{\infty })$ $\forall p\notin B$. Let $K=\bigoplus\limits_{p\in
B}L_{p}.$ In this case, $N/K$ is small in $M/K$ as the set of submodules of $%
\mathbb{Z}(p^{\infty })$ form a chain for all $p\in A$.

Notice that the maximal hollow $\mathbb{Z}$-submodules of $M$ are $\{\mathbb{%
Z}(p^{\infty })\mid p\in A\}$ and they are second, whence $_{\mathbb{Z}}M$
is s-lifting.

Notice that $_{\mathbb{Z}}M$ is second, not semisimple and that $h.\dim (_{%
\mathbb{Z}}M)=\infty .$
\end{ex}

\begin{ex}
\label{Zn[x]}Let $n=p_{1}\cdots p_{n}$ be a product of distinct prime
numbers and consider $M=\mathbb{Z}_{n}[x]$ as a $\mathbb{Z}$-module. Then $M$
is second representable semisimple. Indeed, let $m_{j}=%
\frac{n}{p_{j}}$ for all $j=1,2,...,n$. Set $K_{j_{k}}=\mathbb{Z}m_{j}x^{k}$
for all $j=1, 2, ..., n$ and $k\in \{0, 1, 2, ...\}$. Then $K_{j_{k}}$ is
simple for all $j=1, 2, ..., n$ and $k\in \{0, 1, 2, ...\}$ and $%
K_{j}=\bigoplus\limits_{k=0}^{\infty }K_{j_{k}}$ is $p_{j}$-second for all $%
j=1,2, ..., n$. Hence $M=\bigoplus\limits_{j=1}^{n}K_{j}$ is second
representable while it is semisimple with infinite length.
\end{ex}

The above two examples show also that the finiteness condition on the hollow
dimension in Example \ref{Theorem 4.15} is not necessary.

\begin{ex}
\label{sr-not-ss}Let $n=p_{1}\cdots p_{n}$ be a product of distinct prime
numbers, $p$ any prime number and consider the Abelian group $M=\mathbb{Z}%
_{n}[x]\oplus \mathbb{Z}(p^{\infty })$ as a $\mathbb{Z}$-module. Since $%
\mathbb{Z}_{n}[x]$ is second representable (see Example \ref{Zn[x]}) and $%
\mathbb{Z}(p^{\infty })$ is second, it follows that $_{\mathbb{Z}}M$ is
second representable. Notice that $_{\mathbb{Z}}M$ is neither second nor semisimple.
\end{ex}

As a direct consequence of Lemma \ref{Theorem 1.24} and Example \ref{Theorem
4.15} we obtain:

\begin{cor}
\label{Corollary 4.16}Let $_{R}M$ be s-lifting.

\begin{enumerate}
\item If $M$ has the ACC on supplements (\emph{e.g.} Noetherian), then $M$
is directly second representable.

\item If $M$ has the DCC on supplements (\emph{e.g.} Artinian), then $M$ is
directly second representable.
\end{enumerate}
\end{cor}

\begin{thm}
\label{Theorem 4.18}Let $M$ be an $R$-module.

\begin{enumerate}
\item If $_{R}M$ is finitely generated, supplemented and every maximal
hollow $R$-submodule of $M$ is second, then $M$ is second representable.

\item If $_{R}M$ is coatomic, every maximal $R$-submodule of $M$ has a
supplement in $M,$ every maximal hollow $R$-submodule of $M$ is second and $%
Att^{s}(M)$ is finite, then $M$ is second representable.
\end{enumerate}
\end{thm}

\begin{Beweis}
\begin{enumerate}
\item Since $_{R}M$ is finitely generated and supplemented, $%
_{R}M=\sum\limits_{\lambda \in \Lambda }M_{\lambda }$ a sum of hollow $R$%
-submodules by Lemma \ref{Theorem 1.26} (2). Since $_{R}M$ is finitely
generated, this sum can be taken to be finite and it follows that $M$ is
second representable by Proposition \ref{Proposition 4.9}.

\item Since $_{R}M$ is coatomic and every maximal $R$-submodule of $M$ is
has a supplement in $M,$ it follows by Lemma \ref{Theorem 1.26} (1) that $%
_{R}M=\sum\limits_{\lambda \in \Lambda }M_{\lambda }$ a sum of hollow,
whence second, $R$-submodules of $M.$ Since $Att^{s}(M)$ is finite, it
follows by Proposition \ref{Remark 4.10} that $M$ is second representable.$\blacksquare$
\end{enumerate}
\end{Beweis}

\begin{ex}
\label{Example 4.18}Theorem \ref{Theorem 4.18} provides a non-empty class of
examples of second representable modules. For example, let $p$ be a prime
number and consider the Pr\"{u}fer group $M=\mathbb{Z}(p^{\infty }),$ as a $%
\mathbb{Z}$-module. Notice that
\begin{equation*}  \label{pruf}
\mathbb{Z}({p^{\infty }})=\left\langle g_{1},\text{ }g_{2},\text{ }g_{3},%
\text{ }... \mid g_{1}^{p}=1,\text{ }g_{2}^{p}=g_{1},\text{ }%
g_{3}^{p}=g_{2}, \cdots\right\rangle.
\end{equation*}
Clearly, $M_{\mathbb{Z}}$ is second and supplemented but not finitely generated. This
example shows that the finiteness condition of Theorem \ref{Theorem 4.18}
(1) is not necessary.

Moreover, consider $N=\left\langle g_{k}\right\rangle \leq M$ for some $k\in
\mathbb{N}.$ Observe that $N_{\mathbb{Z}}$ is finitely generated and supplemented and so, by
Theorem \ref{Theorem 4.18}, $N_{\mathbb{Z}}$ is second representable if and only if $N$
is second as it is hollow.
\end{ex}

\begin{defn}
We define a \textit{semisecondary module}%
\index{semisecondary module} as one which is a (possible infinite) sum of
secondary submodules.
\end{defn}

\begin{figure}[tbp]
\centering
\begin{tikzpicture}
\matrix [column sep=7mm, row sep=5mm] {
  \node (se) [draw, shape=rectangle] {Semisecondary}; &
  \node (yw) [draw, shape=circle] {Lifting }; &
  \node (ul) [draw, shape=rectangle] {s-Lifting}; \\
  \node (d1) [draw, shape=circle] {Supplemented}; & &
  \node (d2) [draw, shape=circle] {Semisimple}; \\
  \node (we) [draw, shape=rectangle] {Artinian}; &
  \node (ec) [draw, shape=circle] {Amply supplemented}; &
  \node (pu) [draw, shape=rectangle] {Semisecond}; \\
};
\draw[->, thick] (ec) -- (d1);
\draw[->, thick] (we) -- (ec);
\draw[->, thick] (yw) -- (ec);
\draw[->, thick] (pu) -- (se);
\draw[->, thick] (d2) -- (pu);
\draw[->, thick] (d2) -- (ul);
\draw[->, thick] (ul) -- (yw);
\end{tikzpicture}
\caption{s-lifting position chart}
\label{s-lifting}
\end{figure}
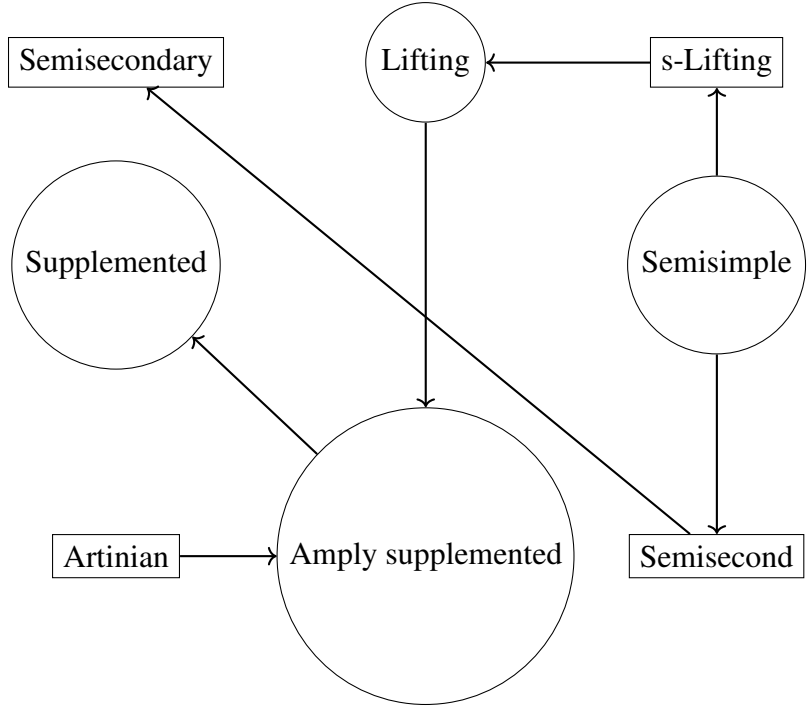

\begin{ex}
Assume that the prime spectrum of $R$ be finite (e.g. $R = \mathds{Z}_n$).
Assume that $M$ is coatomic and amply supplemented over $R$
(e.g. an Artinian module over an Artinian ring) in which the maximal hollow
submodules are second. Then $M$ is second representable by Theorem \ref%
{Theorem 4.18} (2). To show this, let $K\lneq M$ be maximal submodule,
whence there is element $x\in M\backslash K$. So, $K + Rx = M$ as $K$ is
maximal. Since $M_{R}$ is amply supplemented, there is a supplement $N\leq Rx$
of $K$.
\end{ex}

\begin{ex}
\label{Example 4.20}The Abelian group $M=\mathbb{Z}_{12},$ considered as a $%
\mathbb{Z}$-module, has a secondary representation $M=(4)\oplus (3)$ but no
second representation, it has a finite hollow dimension (notice that the
epimorphism
\begin{equation*}
\phi :M\longrightarrow (4)\oplus (6),%
\text{ }x\mapsto 2x
\end{equation*}%
is small and so $h.dim(M)=2$). Observe that $M$ is not s-lifting as the
submodule $(3)$ is maximal hollow but not second. This example shows that
the assumption that $M$ is s-lifting in Theorem \ref{Theorem 4.15} cannot be
dropped.
\end{ex}

\begin{ex}
Let $_{R}M$ have an infinite number of distinct simple $R$-submodules $%
\{S_{1},\cdots \}$ such that $A:=\{Ann(S_{i})\mid i\in \mathbb{N}\}$
is also infinite. The semisimple module $N=\bigoplus\limits_{i=1}^{\infty
}S_{i}$ is not second representable. This example shows that the finiteness
condition on the hollow (uniform) dimension of $_{R}M$ in Theorem \ref%
{Theorem 4.15} cannot be dropped.
\end{ex}

\begin{ex}
A multiplication semisimple module $M=\bigoplus\limits_{i=1}^{\infty }S_{i},$
with infinite number of distinct simple submodules $\{S_{i}\mid i\in \mathbb{%
N}\}$ is not second representable. To prove this, we claim that $%
A:=\{Ann(S_{i})\mid i\in \mathbb{N}\}$ is infinite. Suppose that $%
Ann(S_{i})=P=Ann(S_{j})$ for some $i\neq j$. Since $_{R}M$ is
multiplication, $S_{i}=IM$ for some ideal $I\leq R$, whence $I\subseteq
Ann(S_{j})=P=Ann(S_{i}).$ But this would mean that $S_{i}\nsubseteq IM$ (a
contradiction). Thus, $A$ is infinite as $\{S_{i}\mid i\in \mathbb{N}\}$ is
infinite.
\end{ex}

\section{Second Attached Primes}

Recall that $R$ is a commutative ring. In this section, we investigate the class of (\emph{main}) \emph{second attached primes} of a second representable $R$-module $M$ (with second spectrum $Spec^{s}(M)$).

For every $R$-module $M$, set
\begin{eqnarray*}
Att^{s}(M) &:=&\{Ann(K)\mid K\in Spec^{s}(M)\}
\end{eqnarray*}
If $_{R}M$ is second representable, with a minimal second representation $M=\sum%
\limits_{i=1}^{n}K_{i}$, then the class of main second attached primes of $M$ is given by
\begin{eqnarray*}
att^{s}(M) &=&\{Ann(K_{1}),\cdots ,Ann(K_{n})\}
\end{eqnarray*}

\begin{punto}
\label{second-seconday}Every $p$-second $R$-module is $p$-secondary and
every (minimal) second representation is a (minimal) secondary
representation. So, every second representable $R$-module is secondary
representable and $att^{s}(M)=Att(M)$. A subset $A\subseteq att^{s}(M)$ is
called \emph{isolated}%
\index{isolated subset} iff for any $p\in att^{s}(M)$ with $p\subseteq q$
for some $q\in A$, we have $p\in A$.
\end{punto}

\begin{prop}
\index{$att^{s}(M)$}
\index{$Att^{s}(M)$}
\index{main second attached prime}
\index{second attached prime} \label{Proposition 4.21}Let $_{R}M$ be second
representable, with a minimal second representation $M=\sum%
\limits_{i=1}^{n}K_{i}$, and consider%
\begin{eqnarray*}
att^{s}(M) &:=&\{Ann(K_{1}),\cdots ,Ann(K_{n})\}; \\
Att^{s}(M) &:=&\{Ann(K)\mid K\in Spec^{s}(M)\}
\end{eqnarray*}

\begin{enumerate}
\item $Att^s(M)$ is atomic and $Min(att^s(M)) = Min(Att^s(M))$

\item If there is no small second submodule of $M$, then $Att^{s}(M)$ is
coatomic and $Max(att^{s}(M))=Max(Att^{s}(M))$.
\end{enumerate}
\end{prop}

\begin{Beweis}
\begin{enumerate}
\item \textbf{Claim:} $\bigcap\limits_{p\in att^{s}(M)}p\subseteq q$ for
every $q\in Att^{s}(M)$. Let $a\in \bigcap\limits_{p\in att^{s}(M)}p$. Then $%
a\in Ann(M)$. It is easy to show that $Ann(M)\subseteq \bigcap\limits_{p\in
Att^{s}(M)}p$, hence $$\bigcap\limits_{p\in att^{s}(M)}p=\bigcap\limits_{p\in
Att^{s}(M)}p$$ as $att^{s}(M)\subseteq Att^{s}(M)$. It follows that
$\bigcap\limits_{p\in att^{s}(M)}p\subseteq q$ for all $q\in Att^{s}(M)$.

Now, Suppose that $q\in Min(Att^{s}(M))$. Then $\bigcap\limits_{p\in
att^{s}(M)}p\subseteq q$. Since $att^{s}(M)$ is finite and each element in $%
Att^{s}(M)$ is prime, it follows that $p\leq q$ for some $p\in att^{s}(M)$.
By the minimality of $q$ in $Att^{s}(M)$ and $att^{s}(M)\subseteq Att^{s}(M)$%
, we have $p=q.$ Therefore, $Min(Att^{s}(M))\subseteq Min(att^{s}(M)).$ So, $%
Att^{s}(M)$ is atomic.

For the inverse inclusion, let $p\in Min(att^{s}(M))$. Suppose that $p\notin
Min(Att^{s}(M))$. Then there is $q\in Att^{s}(M)$ such that $q\subsetneq p$.
Since $\bigcap\limits_{p\in att^{s}(M)}p\subseteq q$ and $att^{s}(M)$ is
finite, $p^{\prime }\subseteq q$ for some $p^{\prime }\in att^{s}(M)$, i.e. $%
p^{\prime }\subseteq q\subsetneq p$, which contradicts the minimality of $p$
in $att^{s}(M)$.

\item Assume that there is no small second submodule in $M$.

\textbf{Claim: }For every $p\in Att^{s}(M)$, we have $pM\neq M$ and $%
p\subseteq q$ for some $q\in att^{s}(M)$: Let $p\in Att^{s}(M)$. Then there
is a $p$-second submodule $K\leq M$. Since $K$ is not small in $M,$ there is
a proper submodule $L\lneq M$ such that $K+L=M$ and so $pM=L\neq M$.

Let $p\in Max(Att^{s}(M))$ and assume, without loss of generality that $%
pM=\sum_{i=1}^{m}K_{i}$ with $m\lneq n$ (as $pM\neq M$) and $p\subseteq
Ann(K_{i})$ for all $i\in \{m+1,m+2,\cdots ,n\}$. Since $p\in
Max(Att^{s}(M)) $, $n=m+1$ and $p=Ann(K_{n}),$ i.e. $p\in Max(att^{s}(M))$.
Therefore, $Att^{s}(M)$ is coatomic and $Max(Att^{s}(M))\subseteq
Max(att^{s}(M))$.

For the inverse inclusion, let $q\in Max(att^{s}(M))$. Suppose that $q\notin
Max(Att^{s}(M))$, so that $q\subsetneq p$ for some $p\in Att^{s}(M).$ Then $%
p\subseteq q^{\prime }$ for some $q^{\prime }\in att^{s}(M)$, whence $%
q\subsetneq p\subseteq q^{\prime }$, which contradicts the maximality of $p$
in $att^{s}(M)$. Consequently, $Max(att^{s}(M))\subseteq Max(Att^{s}(M)).$$\blacksquare$
\end{enumerate}
\end{Beweis}

\begin{ex}
\label{Example 4.22}Consider the Abelian group $M=\mathbb{Z}_{n}$ as a $%
\mathbb{Z}$-module. We describe the second spectrum of $M$ and find $%
Att^{s}(M)$ and $att^{s}(M).$
\end{ex}

\begin{ex}
If $n$ is prime, then
\begin{equation*}
Att^{s}(M)=att^{s}(M)=\{(n)\}.
\end{equation*}%
If $n$ is not prime, then consider the prime factorization $%
n=\prod\limits_{i=1}^{k}p_{i}^{n_{i}}$ and let $m_{i}:=n/p_{i}$ for $i\in
\{1,\cdots ,n\}$. Notice that $(m_{i})$ is $p_{i}$-second for all $i\in
\{1,\cdots ,n\}$ and $Att^{s}(M)=\{(p_{1}),\cdots ,(p_{k})\}$.

To find $att^{s}(M)$, we have the following cases:

\textbf{Case 1:} $n_{i}=1$ for all $i\in \{1,\cdots ,n\}$. In this case, $%
M=\sum\limits_{i=1}^{k}(m_{i})$ is a second representation and $%
att^{s}(M)=\{(p_{1}),\cdots ,(p_{k})\}$.

\textbf{Case 2:} $n_{j}>1$ for some $j\in \{1,\cdots ,n\}$. In this case, $%
M$ is not second representable since $\sum\limits_{i=1}^{k}(m_{i})\subseteq
(p_{j})\neq M$.
\end{ex}

\begin{ex}
Let $M$ be a second representable $\mathbb{Z}$-module. Then either $0\in
att^{s}(M)$ or $0\notin att^{s}(M),$ and so by Proposition \ref{Proposition
4.21} we have $Min(Att^{s}(M))=\{0\}$ or $%
Att^{s}(M)=Min(Att^{s}(M))=att^{s}(M)$. In particular, if $M$ is a torsion
module (e.g. $M=\mathbb{Z}_{p}\times \mathbb{Z}_{q}$ for some prime numbers $%
p$ and $q$), then $$Att^{s}(M)=Min(Att^{s}(M))=att^{s}(M).$$
\end{ex}

\begin{thm}
\index{second representation}
\index{directly second representation}
\index{representable module}
\index{$p$-second} \label{Theorem 4.23}Let $_{R}M$ be Noetherian.

\begin{enumerate}
\item Let $p$ be a prime ideal. Then $M$ is $p$-secondary ($p$-second) if
and only if every nonzero submodule of $M$ is $p$-secondary ($p$-second).

\item If $M=\sum\limits_{i=1}^{n}K_{i}$ is a minimal secondary
representation with $K_{i}$ is $p_{i}$-secondary for some prime ideals $%
\{p_{1},\cdots ,p_{n}\}\subseteq Spec(R)$, then $M=\bigoplus%
\limits_{i=1}^{n}K_{i}$.
\end{enumerate}
\end{thm}

\begin{Beweis}
For $_{R}M,$ consider for every $x\in R$ the endomorphism%
\begin{equation*}
a_{M}:M\longrightarrow M,%
\text{ }x\mapsto ax.
\end{equation*}

\begin{enumerate}
\item We prove the result for the case of $p$-secondary modules; the case of
$p$-second modules can be proved similarly.

$(\Rightarrow )$ Let $M$ be a $p$-secondary module for some prime ideal $%
p\leq R$. Let $0\neq K\leq M.$ For any $a\notin p$, we have $aM=M$. Since $%
_{R}M$ is Noetherian, every surjective endomorphism is injective and so $%
a_{M}$ is injective. Hence $a_{M}^{n}$ is injective for any $n$, i.e. $%
a^{n}K\neq 0$ for all $n\in \mathbb{N}$. On the other hand, $aK\subseteq
K=aL $ for some submodule $L$ (as $a_{M}$ is surjective), whence $K\subseteq
L$ and
\begin{equation*}
a_{K}\circ a_{L}:L\longrightarrow aK
\end{equation*}%
is an isomorphism of $R$-modules. Hence $L=K$ and $aK=K$. Therefore $K$
is $p$-secondary. $(\Leftarrow )$ trivial.

\item $M=\sum\limits_{i=1}^{n}K_{i}$ is a minimal secondary representation
with $K_{i}$ is $p_{i}$-secondary for some prime ideals $\{p_{1},\cdots
,p_{n}\}\subseteq Spec(R).$ Let $A=\{1,\cdots ,n\}$.

\textbf{Claim: }For any $j\in A,$ we have $K_{j}\cap \sum\limits_{i\in
A\backslash \{j\}}K_{i}=0.$ Suppose that $K_{j}\cap \sum\limits_{i\in
A\backslash \{j\}}K_{i}\neq 0$ for some $j\in A$. Notice that by (1), $%
N=K_{j}\cap \sum\limits_{i\in A\backslash \{j\}}K_{i}$ is $p_{j}$-secondary.
Set $J:=\{m\in A:p_{m}\nsubseteq p_{j}\}$. For any $m\in J$, there is $%
a_{m}\in p_{m}\backslash p_{j}$. Consider $a=\prod\limits_{m\in J}a_{m}$ and
notice that $a\notin p_{j}$ (as $J$ is finite $a_{m}\in p_{m}\backslash
p_{j} $ for all $m\in J$) and so $aN=N$. Suppose that $\sum\limits_{i\in
A\backslash \{j\}}x_{i}\in N$ such that $x_{i}\in K_{i}$ for all $i\in A$.
Then $a^{l}\sum\limits_{i\in A\backslash \{j\}}x_{i}=\sum\limits_{i\in
A\backslash (\{j\}\cup J)}x_{i}$ for some $l.$ But $a^{l}N=N$ and so $%
N\subseteq \sum\limits_{i\in A\backslash (\{j\}\cup
J)}K_{i}=\sum\limits_{p_{i}\subsetneq p_{j}}K_{i},$ whence $N=K_{j}\cap
\sum\limits_{p_{i}\subsetneq p_{j}}K_{i}$. Since $N\neq 0$, it follows that $%
\{i\in A:p_{i}\subsetneq p_{j}\}\neq \emptyset $. We have the following
cases:

\textbf{Case 1 : }$\{i\in A:p_{i}\subsetneq p_{j}\}=\{h\}$. In this case, $%
N=K_{j}\cap K_{h},$ $N$ is $p_{i}$-secondary and $p_{h}$-secondary at the
same time (a contradiction).

\textbf{Case 2: }$\{i\in A:p_{i}\subsetneq p_{j}\}$\textbf{\ }has more than
one element. In this case, $N=K_{j}\cap \sum\limits_{p_{i}\subsetneq
p_{j}}K_{i}$. Let $p_{h}$ be minimal among all $p_{i}\subseteq p_{j}$, $%
a_{h}\in p_{j}\backslash p_{h}$ and $a_{i}\in p_{i}\backslash p_{h}$ for all
$p_{h}\subsetneq p_{i}\subsetneq p_{j}$ if it exists. Consider $%
b=a_{h}\prod\limits_{p_{h}\subsetneq p_{i}\subsetneq p_{j}}a_{i}$. Then $%
b\in p_{i}$ for all $p_{h}\subsetneq p_{i}\subseteq p_{j}$ and so $%
bK_{h}=K_{h}$. Since $_{R}M$ is Noetherian, $b_{K_{h}}^{t}$ is injective for
every $t$. Hence, for any $\sum\limits_{p_{i}\subsetneq p_{j}}x_{i}\in N$,
we have $b^{t}\sum\limits_{p_{i}\subsetneq p_{j}}x_{i}=b^{k}x_{h}$ for some $%
k$ large enough. But $b_{N}$ is nilpotent as $N$ is $p_{j}$-secondary and $%
b\in p_{j}$. So, $x_{h}=0$ and $N=K_{j}\cap \sum\limits_{p_{i}\subsetneq
p_{j},p_{i}\neq p_{h}}K_{i}$. Also, the set
\begin{equation*}
\{i\in A\mid p_{i}\subsetneq p_{j}\text{ and }p_{i}\neq p_{h}\}
\end{equation*}%
has a minimal element as it is finite. We continue removing the minimal
elements until we arrive at a set containing exactly one element (i.e. Case
1) which yields a contradiction. Therefore, $N=0$.$\blacksquare$
\end{enumerate}
\end{Beweis}

\begin{thm}
\index{second representable} \label{23-loc}Let $_{R}M$ be Noetherian and
Artinian and assume that for any $p\in Max(R)$ the canonical map $\phi
_{p}:M\longrightarrow M_{p}$ is injective. Then $M$ is (second) secondary
representable if and only if $M_{p} $ is (second) secondary representable $%
R_{p}$-module for any $p\in Max(R)$.
\end{thm}

\begin{Beweis}
We prove the result for the case of secondary representation; the case of
second representation can be proved similarly.

Assume that $M_{p}$ is a secondary representable $R_{p}$-module for any
maximal ideal $p\leq R$, say $M_{p}=\sum\limits_{i=1}^{n}K_{i}^{\prime }$ is
a minimal secondary representation for $M_{p}$ where each $K_{i}^{\prime }$
is a secondary submodule of $M_{p}$ and set for all $i\in A=\{1,2,\cdots
,n\}:$
\begin{eqnarray*}
K^{i} &:=&\{x\in M\mid x/s\in K_{i}^{\prime }%
\text{ for some }s\notin p\}; \\
N^{i} &:=&\{x\in M\mid x/1\in K_{i}^{\prime }\}
\end{eqnarray*}%
Then $N^{i}=K^{i}$ and $K_{p}^{i}=K_{i}^{\prime }$. We may write
\begin{equation*}
M_{p}=K_{p}^{1}+K_{p}^{2}+\cdots +K_{p}^{n}.
\end{equation*}%
Consider the canonical map:%
\begin{equation*}
\phi :M\longrightarrow M_{p};x\mapsto x/1.
\end{equation*}%
Then $\phi ^{-1}(K_{i}^{\prime})=K^{i}$ for all $i\in A$ and so $M=$ $%
\sum\limits_{i=1}^{n}K^{i}$.

\textbf{Claim:} for any $a\in R$ and all $i\in A$, the map%
\begin{equation*}
\phi ^{a,i}:K^{i}\longrightarrow K^{i};x\mapsto ax
\end{equation*}

is nilpotent or surjective. To show this, suppose that $\phi ^{a,i}$ is not
surjective. Then $\phi ^{a,i}$ is not injective since $K^{i}$ is Artinian
(every injective endomorphism of an Artinian module is surjective). Since $%
\phi ^{a,i}$ is not injective, \cite[Proposition 3.9]{AM1969} implies that $%
\phi _{p}^{a,i}$ is not injective for some maximal ideal $p$. So, $\phi
_{p}^{a,i}$ is not surjective (any surjective endomorphism of a Noetherian
module is injective). It follows that $\phi _{p}^{a,i}$ is nilpotent, i.e.
for some $n$ we have $a^{n}x/1=0$ for all $x\in K^{i}$ whence $a^{n}x=0$ for
all $x\in K^{i}$ by our assumption.

The converse is clear (see Remark \ref{Remark 4.6}(5) for the second
representation case).
$\blacksquare$
\end{Beweis}

\begin{ex}
\label{Example 5.18}Consider the $\mathbb{Z}$-module $M=\prod\limits_{i=1}^{%
\infty }\mathbb{Z}_{p_{i}p_{i}^{\prime }}$, where $p_{i}$
and $p_{i}^{\prime }$ are primes and $p_{i}\neq p_{j}$, $p_{i}^{\prime }\neq
p_{j}^{\prime }$ for all $i\neq j\in \mathbb{N}$ and $p_{i}^{\prime }\neq
p_{j}$ for any $i$ and $j$. Let the simple $\mathbb{Z}$-modules $K_{p_{i}}$
and $K_{p_{i}^{\prime }}$ be such that $(0:K_{p_{i}})=(p_{i})$ and $%
(0:K_{p_{i}^{\prime }})=(p_{i}^{\prime })$, so
\begin{equation*}
M=\bigoplus\limits_{i=1}^{\infty }K_{p_{i}}\oplus \bigoplus_{i=1}^{\infty
}K_{p_{i}^{\prime }}.
\end{equation*}%
Every second $\mathbb{Z}$-submodule of $M$ is simple, while $_{\mathbb{Z}}M$
is not multiplication.
\end{ex}

\begin{ex}
\label{Example 4.24}

\begin{enumerate}
\item There exists an $R$-module $M$ which is $p$-secondary but not
Noetherian, while every submodule of $M$ is $p$-secondary. Appropriate
semisimple modules with infinite lengths provide a source for such modules,
see Example \ref{Example 5.18}.

\item If $_{R}M$ is a Noetherian or an Artinian $R$-module with no zero
divisors, then $M$ satisfies the conditions of Theorem \ref{23-loc}.
\end{enumerate}
\end{ex}

\begin{thm}
\index{multiplication module}
\index{semisecond}
\index{$Att^{s}(M)$}
\index{$att^{s}(M)$} \label{Theorem 4.25}

\begin{enumerate}
\item A multiplication $R$-module $M$ is semisecond (resp. second
representable) if and only if each nonzero submodule of $M$ is semisecond
(resp. second representable).

\item A multiplication $R$-module $M$ which is not hollow is semisecond
(resp. second representable) if and only if each non-small proper submodule
of $M$ is semisecond (resp. second representable).

\item The following conditions are equivalent for a second representable $R$%
-module $M$ with a minimal second representation $\sum\limits_{i=1}^{n}K_{i}$%
.

\begin{enumerate}
\item $M$ is multiplication.

\item $Att^{s}(M)=att^{s}(M)=Min(att^{s}(M))$ and every nonzero submodule of
$M$ has a second representable $\sum\limits_{j\in A}K_{j}$ for some $%
A\subseteq \{1,\cdots ,n\}$.
\end{enumerate}

\item If $M$ is semisimple second representable, then the following
conditions are equivalent:

\begin{enumerate}
\item $M$ is multiplication.

\item The elements of $att^{s}(M)$ are incomparable and every second
submodule of $M$ is simple.
\end{enumerate}

\item The following conditions are equivalent for an atomic module $_{R}M:$

\begin{enumerate}
\item $M$ is semisimple.

\item $M=\sum\limits_{i\in A}K_{i}$ where every submodule of $K_{i}$ is $%
p_{i}$-second for some prime ideal $p_{i},$ $i\in A$.
\end{enumerate}
\end{enumerate}
\end{thm}

\begin{Beweis}
\begin{enumerate}
\item Let $_{R}M$ be multiplication. Suppose that $M=\sum\limits_{i\in A}K_{i}$ is a
semisecond representation of $M$. Let $0\neq K\leq M$, whence $K=IM$ for
some ideal $I\leq R$. Suppose that $I\nleq p_{j}$ for all $j\in B\subseteq A$
and $I\leq p_{i}$ for all $i\in A\backslash B$. Then $IM=\sum\limits_{i\in
B}K_{j}$ and so $K$ is semisecond.

\item Assume that $_{R}M$ is multiplication and not hollow. Then $%
M=K_{1}+K_{2}$, where $K_{1}$ and $K_{2}$ are proper non-small submodules
of $M$. Notice that $K_{1}$ and $K_{2}$ are semisecond, whence $M$ is semisecond.

\item Let $M=\ \sum\limits_{i=1}^{n}K_{i}$ be a minimal second
representation.

$(a\Rightarrow b)$ Assume that $M$ is a multiplication module. Let $p\in
Att^{s}(M)$, whence there is $K\leq M$ which is $p$-second. Since $_{R}M$ is
multiplication, $K=IM$ for some ideal $I\leq R$. Assume, without loss of generality,
that $IM=\sum\limits_{i=1}^{m}K_{i}$ where $I\subseteq p_{i}$ for all
$i\in \{m+1,\cdots ,n\}$. Since $K$ is $p$-second, it follows by
Proposition \ref{Proposition 4.2} that $m=1$ and $p=p_{i}$ for some $i\in
\{1,\cdots ,n\}$ and so $Att^{s}(M)=att^{s}(M)$. If $K\neq 0$, then there
is an ideal $J\leq R$ such that $K=JM=J\sum\limits_{i=1}^{n}K_{i}=\sum%
\limits_{i\in A}K_{i}$, where
\begin{equation*}
A=\{i\in \{1,2,\cdots ,n\}\mid J\nsubseteq p_{i}\}.
\end{equation*}
Assume that $p_{i}\subsetneq p_{j}$ and $K_{j}=IM$, but $IM=I\sum%
\limits_{i=1}^{n}K_{i}$. By the minimality of $\sum\limits_{i=1}^{n}K_{i}$
and using Proposition \ref{Proposition 4.2}, $I\subseteq p_{i}$ and so $%
I\subseteq p_{j}$ (a contradiction).

$(b\Rightarrow a)$ Assume that $K\leq M$ is a nonzero submodule. By our
assumption, $K=\sum\limits_{j\in A}K_{j}$ with some $A\subseteq \{1,\cdots
,n\}$. Assume, without loss of generality, that $K=\sum\limits_{i=1}^{m}K_{i}$%
. Let $I=\bigcap\limits_{i=m+1}^{n}p_{i}$. Since $I\subseteq p_{i}$ for all $%
i\in \{m+1,\cdots,n\}$ and since every element in $Att^{s}(M)$ is minimal
and strongly irreducible and $Att^{s}(M)$ is finite, it follows that $%
I\nsubseteq p_{i}$ for all $i=1,\cdots ,m$. Therefore, $K=IM$.
Consequently, $_{R}M$ is multiplication.

\item Let $_{R}M$ be semisimple with a second presentation $%
M=\sum\limits_{i=1}^{n}K_{i}.$ It follows that $Att^{s}(M)=att^{s}(M)$.
Apply now (3) and observe that every nonzero submodule of $M$ is the second
representable module $\sum\limits_{j\in A}K_{j}$ for some $A\subseteq
\{1,\cdots ,n\}$ if and only if every second submodule of $M$ is simple.

\item Assume that $M$ is atomic and semisecond, say $M=\sum\limits_{i\in
A}K_{i}$ where every submodule of $K_{i}$ is $p_{i}$-second for some prime
ideal $p_{i},$ $i\in A$. For each $i\in \{1,\cdots ,n\}$ we set
\begin{equation*}
H_{i}=\{S:S \text{ is simple in }K_{i}\}.
\end{equation*}
\textbf{Claim:} $K_{i}=\bigoplus\limits_{S\in H_{i}}S$. If not, then there
exits $x\in K_{i}\backslash \bigoplus_{S\in H_{i}}S$ and so $Rx$ must
contain some $S\in H_{i}$ and so there is $a\in R$ such that $ax\in S$,
whence $S=Rax$. But $Rx$ is $p_{i}$-second and $Rax=S\neq 0$, whence $%
S=Rax=Rx$ and so $x\in S$ (a contradiction). It follows that $M$ is a sum of
simple submodules. The converse is trivial.$\blacksquare$
\end{enumerate}
\end{Beweis}

\begin{ex}
\index{atomic} Every semisecond atomic Noetherian $R$-module is semisimple.
This follows directly from Theorem \ref{Theorem 4.23} (1) and Theorem \ref%
{Theorem 4.25} (5).
\end{ex}

\begin{thm}
\index{injective module}
\index{main second attached prime}
\index{$att^{s}(M)$} \label{Theorem 4.26}Let $E$ be an injective $R$-module
and $\bigcap\limits_{i=1}^{n}p_{i}=0$ where $p_{1},\cdots ,p_{n}$ are
incomparable prime ideals.

\begin{enumerate}
\item $E$ is second representable and $att^{s}(E)\subseteq
\{p_{1},p_{2},\cdots ,p_{n}\}$.

\item $Ann(E)=0$ if and only if $att^{s}(E)=\{p_{1},\cdots ,p_{n}\}$.
\end{enumerate}
\end{thm}

\begin{Beweis}
\begin{enumerate}
\item Set $E[I]:=(0:_{E}I)$. By \cite[Lemma 2.1]{S1976}, $E[p_{i}]=0$ or $%
E[p_{i}]$ is $p_{i}$-secondary for all $i\in \{1,2,\cdots ,n\}$. Assume that
$E[p_{i}]\neq 0$.

\textbf{Claim:} $E[p_{i}]$ is second. Let $I\leq R$. If $I\subseteq p_{i}$,
then $IE[p_{i}]=0$. If $I\nsubseteq p_{i}$, then there is $a\notin p_{i}$.
Since $E[p_{i}]$ is $p_{i}$-secondary, $aE[p_{i}]=E[p_{i}]$ and so $%
IE[p_{i}]=E[p_{i}]$.

Notice that $E=\sum\limits_{i=1}^{n}E[p_{i}]$ follows from \cite[Lemma 2.2]%
{S1976} as $E=E[0]=E[\bigcap\limits_{i=1}^{n}p_{i}]=\sum%
\limits_{i=1}^{n}E[p_{i}]$. Hence, $E$ is second representable and $%
att^{s}(E)\subseteq \{p_{1},\cdots ,p_{n}\}$.

\item $(\Rightarrow )$ Assume that $Ann(E)=0$. We want to show that $%
E=\sum\limits_{i=1}^{n}E[p_{i}]$ is a minimal second representation. By \cite%
[Lemma 2.2]{S1976}, $E[\bigcap\limits_{i\neq j}p_{i}]=\sum\limits_{i\neq
j}E[p_{i}]$ for any $j\in \{1,\cdots ,n\}$ and so $E\neq
\sum\limits_{i\neq j}E[p_{i}]$ for any $j\in \{1,\cdots ,n\}$. Otherwise, $%
E=E[\bigcap\limits_{i\neq j}p_{i}]=\sum\limits_{i\neq j}E[p_{i}]$ for some $%
j\in \{1,\cdots ,n\}$ and so $\bigcap\limits_{i\neq j}p_{i}$ annihilates
every element in $E$ for some $j\in \{1,\cdots ,n\}$. But $Ann(E)=0$,
whence $\bigcap_{i\neq j}p_{i}=0$, which contradicts the fact that $%
p_{i}\nsubseteq p_{j}$ for all $i\in \{1,\cdots ,n\}\backslash \{j\}$.
Therefore, $att^{s}(E)=\{p_{1},\cdots ,p_{n}\}$.

$(\Leftarrow )$ Assume that $att^{s}(E)=\{p_{1},\cdots ,p_{n}\}$.
Suppose that $0\neq a\in Ann(E).$ Then $a\notin p_{j}$ for some $j\in
\{1,\cdots ,n\}$ and so $E[p_{j}]=0$ or $aE[p_{j}]=E[p_{j}].$ Since $p_{j}$
is a second attached prime, $E[p_{j}]\neq 0$. Therefore, $%
aE=\sum_{i=1}^{n}aE[p_{i}]\neq 0$.$\blacksquare$
\end{enumerate}
\end{Beweis}

\begin{ex}
Injective modules over Artinian rings are second representable. So, any
module over an Artinian ring is embedded in a second representable one,
namely, the injective hull of it.
\end{ex}

\addcontentsline{toc}{section}{\protect\numberline{}{Index}} \printindex


\begin{thebibliography}{99}

\bibitem{Abu} J. Abuhlail, \emph{Zariski topologies for coprime and second
submodules}, Algebra Colloq., 22, 47-72 (2015).

\bibitem{AK2012} A. Altman and S. Kleiman, \emph{A Term of Commutative
Algebra}, MIT (2012).

\bibitem{Ann2002} S. Annin, \emph{Associated and Attached Primes over
Noncommutative Rings,} Ph.D. Dissertation, University of California at
Berkeley (2002).

\bibitem{AM1969} M. Atiyah and I. Macdonald, \emph{Introduction to
Commutative Algebra}, Addison-Wesley Publishing Co., Reading,
Mass.-London-Don Mills, Ont. (1969).

\bibitem{BH1990} Y. Baba and M. Harada, \emph{On Almost $M$-Projectives and Almost $M$-Injectives}, Tsukuba J. Math. 14, 53-69 (1990).

\bibitem{B2009} M. Baig, \emph{Primary Decomposition and Secondary
Representation of Modules over a Commutative Ring}, Thesis, Georgia State
University (2009).

\bibitem{JCNR} J. Clark, Ch. Lomp, N. Vanaja and R. Wisbauer, \emph{Lifting
Modules, Supplements and Projectivity in Module Theory, }(Frontiers in
Mathematics), Birkh\"{a}user (2006).

\bibitem{PRU} Laszlo Fuchs. \emph{Infinite abelian groups}. Academic Press,
INC. New York 1973.

\bibitem{C1981} L. Chambles, \emph{Coprimary decomposition, N-dimension and
divisibility: Application to Artinian modules, }Commun. Algebra 9, 1131-1146
(1981).

\bibitem{K1973} D. Kirby, \emph{Coprimary decomposition of Artinian modules}%
, J. London Math. Soc. 6, 571-576 (1973).

\bibitem{M1973} I. G. Macdonald, \emph{Secondary representation of modules
over a commutative ring}, Symp. Math. 11, 23-43 (1973).

\bibitem{S1976} R. Y. Sharp, \emph{Secondary representations for injective
modules over commutative Noetherian rings}, Proc. Edinb. Math. Soc. 20,
143-151 (1976).

\bibitem{ST2009} Semra Tekin, \emph{Modules with Coprimary Decompositions},
Izmir Institute of Technology, M.Sc. Thesis (2009).

\bibitem{Y2001} S. Yassemi, \emph{The dual notion of prime submodules},
Arch. Math. (Brno) 37, 273-278 (2001).

\bibitem{Y1995} S. Yassemi, \emph{Coassociated primes}, Commun. Algebra 23,
1473-1498 (1995).

\bibitem{Wij2006} I. Wijayanti, \emph{Coprime modules and comodules}, Ph.D.
Dissertation, Heinrich-Heine Universit\"{a}t, D\"{u}sseldorf (2006).

\bibitem{W1991} R. Wisbauer, \emph{Foundations of Module and Ring Theory, A
Handbook for Study and Research}, CRC Press, (1991).
\end{thebibliography}
\end{document}